\documentclass[12pt]{article}
\catcode`\@=11

\newtheorem{lemma}{Lemma}
\newtheorem{proposition}{Proposition}

\newtheorem{corollary}{Corollary}
\def\demo{\noindent{\bf Proof .-}}
\def\section{\@startsection {section}{1}{\z@}{-3.5ex plus -1ex
minus-.2ex}{2.3ex plus .2ex}{\normalsize\bf}}
\def\codim{{\rm codim}\,}

\def\bz{\hbox{\it Z\hskip -4pt Z}}

\newcommand{\het}{H_{\rm et}}
\newcommand{\hc}{H_{\rm c}}
\begin{document}
\begin{center}
{\Large\bf \textsc{On a special class of simplicial toric varieties}}
\end{center}
\vskip.5truecm
\begin{center}
{Margherita Barile\footnote{Partially supported by the Italian Ministry of Education, University and Research.}\\ Dipartimento di Matematica\\ Universit\`{a} di Bari\\ Via E. Orabona 4\\70125 Bari, Italy}
\end{center}
\vskip1truecm
\noindent
{\footnotesize{\bf Abstract}  We show that for all $n\geq 3$ and all primes $p$ there are infinitely many simplicial toric varieties of codimension  $n$ in the $2n$-dimensional affine space whose minimum number of defining equations is equal to $n$ in characteristic $p$, and lies  between $2n-2$ and $2n$ in all other characteristics. In particular, these are new examples of varieties which are set-theoretic complete intersections only in one positive characteristic.\newline
 Moreover, we show that the minimum number of binomial equations which define these varieties in all characteristics is 4 for $n=3$ and $2n-2+{n-2\choose 2}$ whenever $n\geq 4$. \vskip0.5truecm
\noindent
}

\section*{Introduction}
Let $K$ be an algebraically closed field, and let $V$ be an affine variety in $K^N$. The arithmetical rank (ara) of $V$ is defined as the least number of equations that are needed to define $V$ set-theoretically as a subvariety of $K^N$. In general we have that ara\,$V\geq\codim\,V$. If equality holds, $V$ is called a {\it set-theoretic complete intersection}. In general, the arithmetical rank of a variety may depend upon the characteristic of the ground field: but not many examples of this kind are known so far. The first ones to be found were the determinantal varieties of a symmetric matrix considered in \cite{B0}. The paper \cite{BL} presents an infinite class of simplicial toric varieties of codimension 2 which are set-theoretic complete intersection only in one positive characteristic. The same property has been shown in  \cite{B2} for the Veronese varieties whose degree is  a prime power. These have arbitrarily high codimensions, but only the codimensions of a special form are represented; moreover, in each possible codimension there are only finitely many examples. \newline
In the present paper we will show that for every prime $p$ and every  codimension $n\geq3$, there are infinitely many examples of simplicial toric varieties in $K^{2n}$ which are set-theoretic complete intersections  (i.e. set-theoretically defined by $n$ equations) only in characteristic $p$. This result completes the one in \cite{BL}, but with an interesting difference: as it was proven in \cite{BMT}, the simplicial toric varieties of codimension 2 are always (regardless of their dimension) set-theoretically defined by 3 equations, whereas the arithmetical rank of the varieties that we are going to introduce here dramatically increases from $n$ to a number lying between $2n-2$ and $2n$ when the characteristic is different from $p$.\newline
It is well-known that the defining ideal of every toric variety is generated by binomials. According to a definition introduced by Thoma in \cite{Th}, the {\it binomial arithmetical rank} (bar) of a toric variety $V$ is the least number of binomial equations which are needed to define $V$. Obviously ara\,$V\leq\,{\rm bar}\,V$. A complete characterization of the cases where equality holds is given in \cite{BMT2}. For the varieties $V$ that will be presented in this paper, we show that equality holds only for the single positive characteristic where $V$ is a set-theoretic complete intersection; in this case bar\,$V=\codim V$, and $V$ is therefore called a {\it binomial set-theoretic complete intersection}.  In the remaining characteristics, ara\,$V$ and bar\,$V$ differ in a substantial way, since  we have that bar\,$V=4$ if $n=3$, and bar$\,V=2n-2+{n-2\choose 2}$ for all $n\geq 4$. Thus we have new examples of affine toric varieties of any codimension greater than or equal to 3 such that, in all but one positive characteristic,  the minimum number of defining equations cannot be attained by systems of binomial equations. This property is known to be true, in characteristic zero, for certain projective toric curves in ${\bf P}_K^3$, among which the famous Macaulay's curve (see \cite{Th0}); other special classes of projective toric varieties fulfilling this property in all characteristics have been recently presented in \cite{KMT}, Section 5, and in \cite{K}, Sections 5--6.

\section{Preliminaries} 
A {\it monomial} in  a polynomial ring is a product  of indeterminates. Given a monomial $M$, we define the {\it support} of $M$, denoted supp\,$(M)$, as the set of indeterminates which divide $M$.\newline
A {\it binomial} is the difference of two distinct monomials $M, M'$: these will be called  the monomials of the binomial $B=M-M'$.  This binomial is called {\it monic} in the indeterminate $z$ if supp\,$(M)=\{z\}$.
\par
Let $n\geq3$ be an integer. Moreover, let $b_1,\dots, b_{n-2}$ be nonnegative integers,  $p$ be a prime, and $\ell, a,d$,  $c_1,\dots, c_{n-2}$ be positive integers such that
\begin{list}{}{}
\item{(I)} $p$ does not divide any of the $c_i$;
\item{(II)}  $a$ and $d$ are coprime;
\item{(III)} there are positive integers $g, h$ such that $p^{\ell}=ag+dh$.
\end{list}
\par\medskip\noindent
 Consider the affine {\it simplicial toric} variety $V\subset K^{2n}$ admitting the following parametrization:
$$V:\left\{
\begin{array}{rcl}
x_1&=&u_1\\
&\vdots&\\
x_{n-2}&=&u_{n-2}\\
x_{n-1}&=&u_{n-1}^{p^{\ell}}\\
x_n&=&u_n^a\\
y_1&=&u_1^{b_1}u_{n-1}^{c_1}\\
&\vdots&\\
y_{n-2}&=&u_{n-2}^{b_{n-2}}u_{n-1}^{c_{n-2}}\\
y_{n-1}&=&u_n^d\\
y_n&=&u_{n-1}u_n\\
\end{array}\right..$$
\noindent
We have that $\codim V=n$. Let $I(V)$ be the defining ideal of $V$ in the polynomial ring $R=K[x_1,\dots, x_n, y_1,\dots, y_n]$. Ideal $I(V)$ is prime and is generated by binomials.
The binomials in $I(V)$ are those of the form:  
$$
x_1^{\alpha_1^+}\cdots x_n^{\alpha_n^+}y_1^{\beta_1^+}\cdots y_n^{\beta_n^+}-
x_1^{\alpha_1^-}\cdots x_n^{\alpha_n^-}y_1^{\beta_1^-}\cdots y_n^{\beta_n^-}$$
where $({\alpha_1,\dots,\alpha_n,\beta_1,\dots,\beta_n})\in\bz^{2n}\setminus\{\bf{0}\}$  is such that
\begin{eqnarray*}&&\alpha_1{\bf e}_1+\cdots+\alpha_{n-1}p^{\ell}{\bf e}_{n-1}+\alpha_na{\bf e}_n+\\
&&\beta_1(b_1{\bf e}_1+c_1{\bf e}_{n-1})+\cdots+\beta_{n-2}(b_{n-2}{\bf e}_{n-2}+c_{n-2}{\bf e}_{n-1})+\\
&&\beta_{n-1}d{\bf e}_n+\beta_n({\bf e}_{n-1}+{\bf e}_n)={\bf0}
\end{eqnarray*}
\noindent
where ${\bf e}_1,\dots,{\bf e}_n$ is the standard basis of $\bz^n$, and we have set $\alpha_i^+=\max\{\alpha_i, 0\}$, $\alpha_i^-=\max\{-\alpha_i, 0\}$, $\beta_i^+=\max\{\beta_i, 0\}$ and  $\beta_i^-=\max\{-\beta_i, 0\}$.\par
The next result will be useful in one of the proofs of Section 3.
\begin{lemma}\label{preliminary1} Let $B$ be a binomial of $I(V)$, and let $M, M'$ be its monomials. Then the following conditions hold.
\begin{list}{}{}
\item{(i)} For all $i=1,\dots, n-2$, $x_i$ divides $M$ if and only if $y_i$ divides $M'$.
\item{(ii)} If one of the indeterminates $x_n, y_{n-1}, y_n$ divides $M$, then one of the remaining two divides $M'$.
\item{(iii)} If $x_{n-1}$ divides $M$, then one of $y_1,\dots, y_{n-2}, y_n$ divides $M'$. 
\end{list}
\end{lemma}
\noindent
\demo We prove (i). Let $i$ be any index with $1\leq i\leq n-2$. Let $u_i=0$ and $u_k=1$ for $k\ne i$. Let ${\bf x}=(\bar x_1,\dots, \bar x_n, \bar y_1,\dots, \bar y_n)\in V$ be the  point corresponding to this choice of parameters. Then 
\begin{equation}\label{ast}\bar x_i=\bar y_i=0\qquad\mbox{ and }\qquad\bar x_k=\bar y_k=1\qquad\mbox{ for all }k\ne i.\end{equation}
\noindent
 Since $B\in I(V)$, we have that $B({\bf x})=0$. Suppose that $x_i$ divides $M$. Then $M({\bf x})=0$, so that $M'({\bf x})=0$. On the other hand, by irreducibility, $x_i$ does not divide $M'$. In view of (\ref{ast}) it follows that $y_i$ divides $M'$. The proof of the converse is identical. \par
Claims (ii) and (iii) can be shown by similar arguments, by selecting certain points of $V$ and using the fact that monomial $B$ vanishes at these points. For the proof of (ii) take $u_n=0$ and $u_k=1$ for all $k\ne n$, for the proof of (iii) take $u_{n-1}=0$ and $u_k=1$ for all $k\ne n-1$. 
\par\smallskip\noindent
Consider the following binomials of $R$:
\begin{eqnarray}
F_i&=&y_i^{p^{\ell}}-x_i^{p^{\ell} b_i}x_{n-1}^{c_i}\qquad\qquad(i=1,\dots, n-2)\label{Fi}\\
F_{n-1}&=&y_{n-1}^a-x_n^d\label{Fn-1}\\
F_n&=&y_n^{p^{\ell}}-x_{n-1}x_n^{g}y_{n-1}^{h}\label{Fn}
\end{eqnarray}
\noindent
In view of (III) it easily follows  that $F_1,\dots, F_n\in I(V)$. \newline
The next result can be shown by the same arguments as Lemma \ref{preliminary1}. Nevertheless, we give its proof for the sake of completeness.  
\begin{lemma}\label{preliminary2} For all $i=1,\dots, n-1$, $F_i$ is the only irreducible binomial in $I(V)$ that is monic in $y_i$. 
\end{lemma} 
\demo Consider an index $i$ with $1\leq i\leq n-2$. Let $u_i=u_{n-1}=1$ and $u_k=0$ for all $k\ne i, n-1$, and let ${\bf x}=(\bar x_1,\dots, \bar x_n, \bar y_1,\dots, \bar y_n)$ be the corresponding point of $V$. Then
\begin{equation}\label{xy} \bar x_k=0\quad\mbox{ for all }k\ne i, n-1\qquad\mbox{ and }\quad\bar y_k=0\quad\mbox{ for all }k\ne i.\end{equation}
\noindent
Let $B=M-M'\in I(V)$, where $M$ and $M'$ are monomials, and suppose that $M=y_i^{\alpha_i}$, with $\alpha_i>0$. Then $M({\bf x})=1$, so that, being $B({\bf x})=0$, we have $M'({\bf x})=1$. In view of (\ref{xy}) and irreducibility, this implies that supp\,$(M')\subset\{x_i, x_{n-1}\}$. So let $M'=x_i^{\beta_i}x_{n-1}^{\gamma_i}$. Then 
\begin{equation}\label{beta} \alpha_ib_i=\beta_i,\end{equation}
\begin{equation}\label{gamma} \alpha_ic_i=\gamma_ip^{\ell}.\end{equation}
\noindent
From (\ref{gamma}) and (I) we deduce that $p^{\ell}$ divides $\alpha_i$; let $\alpha'_i=\frac{\alpha_i}{p^{\ell}}$. Then, by (\ref{gamma}), $\alpha'_ip^{\ell}c_i=\gamma_ip^{\ell}$, whence $\alpha'_ic_i=\gamma_i$. Thus, by (\ref{beta}),
$$B=y_i^{\alpha'_ip^{\ell}}-x_i^{\alpha'_ip^{\ell}b_i}x_{n-1}^{\alpha'_ic_i}.$$
\noindent
But irreducibility implies that $\alpha'_i=1$. Hence $B=F_i$, as required. The proof for $i=n-1$ is analogous: it suffices to consider the point of $V$ corresponding to $u_n=1$ and $u_k=0$ for $k\ne n$ and to argue as above. 
\par\medskip\noindent

\section{The defining equations}
In this section we will explicitly give the binomial defining equations for variety $V$. We will have to distinguish two cases, according to the characteristic of the ground field $K$. 
\begin{proposition}\label{equationsp} Suppose that {\rm char}\,$K=p$. Then $V$ is set-theoretically defined by
$$F_1=\cdots=F_n=0.$$
\end{proposition}
\demo We only have to prove that every ${\bf x}\in K^{2n}$ fulfilling the given equations belongs to $V$. So let ${\bf x}=(\bar x_1,\dots, \bar x_n, \bar y_1,\dots, \bar y_n)\in K^{2n}$ be such that $F_i({\bf x})=0$ for all $i=1,\dots, n$. Set
\begin{equation}\label{xi}u_i=\bar x_i\qquad\mbox{ for }i=1\dots, n-2,\end{equation}
\noindent and, moreover, let $u_{n-1}, v_n\in K$ be such that 
\begin{equation}\label{xn-1}\bar x_{n-1}=u_{n-1}^{p^{\ell}}\end{equation}
\noindent
and $\bar x_n=v_n^a$. We show that up to replacing $v_n$ with another $a$-th root $u_n$ of $\bar x_n$ in $K$,  ${\bf x}$ fulfills the parametrization of $V$ given above. 
Condition $F_{n-1}({\bf x})=0$ implies that $\bar y_{n-1}^a=v_n^{ad}$, i.e., $\bar y_{n-1}=v_n^d\omega$ for some $\omega\in K$ such that $\omega^a=1$. By virtue of (II) there are integers $r, s$ such that $1=ar+ds$. Set $\eta=\omega^s$. Then $\eta^a=1$ and $\eta^d=\omega^{ds}=\omega^{ar+ds}=\omega.$ Put $u_n=v_n\eta$. Then 
\begin{equation}\label{xn}\bar x_n=v_n^a=v_n^a\eta^a=u_n^a, 
\end{equation}
\noindent
and 
\begin{equation}\label{yn-1}\bar y_{n-1}=v_n^d\omega=(v_n\eta)^d=u_n^d, 
\end{equation}
\noindent
Furthermore, for all $i=1,\dots, n-2$, $F_i({\bf x})=0$ implies that $\bar y_i^{p^{\ell}}=u_i^{p^{\ell}b_i}u_{n-1}^{p^{\ell}c_i}$, i.e., 
\begin{equation}\label{yi}\bar y_i=u_i^{b_i}u_{n-1}^{c_i}.\end{equation}
\noindent
Finally, in view of (III), $F_n({\bf x})=0$ implies that $\bar y_n^{p^{\ell}}=u_{n-1}^{p^{\ell}}u_n^{ag}u_n^{dh}=u_{n-1}^{p^{\ell}}u_n^{p^{\ell}}, $ i.e., 
\begin{equation}\label{yn} \bar y_n=u_{n-1}u_n.\end{equation}
\noindent
From (\ref{xi})--(\ref{yn}) it follows that ${\bf x}$ fulfills the required parametrization. This completes the proof.
\par\medskip\noindent
We have thus proven that $V$ is set-theoretically defined by $n$  binomial equations, i.e., we have the following
\begin{corollary}\label{corollary1}
If char\,$K=p$,  variety $V$ is a binomial set-theoretic complete intersection.
\end{corollary}
\par\medskip\noindent
We will show that the above corollary does not extend to the characteristics different from $p$. In general $F_1,\dots, F_n$ do not suffice to define $V$ set-theoretically:  more binomial equations are needed. We are going to define these first. 
By virtue of (I), for all indices $i,j$ with $1\leq i<j\leq n-2$, there are positive integers $g_{ij}$, not divisible by $p$,  such that 
\begin{equation}\label{gi}  c_ig_{ij}+c_jg_{ji}=r_{ij}p^{m_{ij}}\end{equation}
\noindent
for suitable positive integers $r_{ij}, m_{ij}$. Note that  $g_{ij}$ and $g_{ji}$ can be chosen large enough so as to have $m_{ij}\geq \ell$.
By (I), for all indices $i$, $1\leq i\leq n-2$, there are also positive integers $h_i, k_i$, not divisible by $p$, such that 
\begin{equation}\label{hiki} c_ih_i+k_i=r_ip^{m_i}.
 \end{equation}
\noindent
for suitable positive integers $r_i, m_i$. Note that  $h_i$ and $k_i$ can be chosen large enough so as to have $m_i\geq \ell$. By (II) there are integers $s_i, t_i$ such that 
\begin{equation}\label{siti} as_i+dt_i=k_i.
\end{equation}
\noindent  
Up to replacing $h_i, k_i$ with  larger numbers we may assume that $s_i, t_i$ are nonnegative. Put
$$G_{ij}=y_i^{g_{ij}}y_j^{g_{ji}}-x_i^{b_ig_{ij}}x_j^{b_jg_{ji}}x_{n-1}^{r_{ij}p^{m_{ij}-\ell}}\qquad(1\leq i<j\leq n-2),$$
\noindent 
and
$$H_i=y_i^{h_i}y_n^{k_i}-x_i^{b_ih_i}x_{n-1}^{r_ip^{m_i-\ell}}x_n^{s_i}y_{n-1}^{t_i}\qquad(i=1,\dots, n-2).$$
\noindent
From (\ref{gi}), (\ref{hiki}) and (\ref{siti}) it follows that $G_{12}, G_{13},\dots, G_{n-3\,n-2}, H_1,\dots, H_{n-2}\in I(V)$.
\begin{proposition}\label{equationsnonp} Variety $V$ is set-theoretically defined by 
$$F_1=F_2=F_3=H_1=0\qquad\qquad\mbox{if }n=3,$$
\noindent and by 
$$F_1=\cdots=F_n=G_{12}=\cdots=G_{n-3\,n-2}=H_1=\cdots=H_{n-2}=0\qquad\mbox{if }n\geq4.$$
\end{proposition}
\demo Suppose that ${\bf x}=(\bar x_1,\dots, \bar x_n, \bar y_1,\dots, \bar y_n)\in K^{2n}$ fulfills the system of equations given in the claim. Set $u_i=\bar x_i$ for all $i=1,\dots, n-2$. As in the proof of Proposition \ref{equationsp}, $F_{n-1}({\bf x})=0$ implies that there is $u_n\in K$ such that $\bar x_n=u_n^a$ and $\bar y_{n-1}=u_n^d$.  Let $v_{n-1}\in K$ be such that $\bar x_{n-1}=v_{n-1}^{p^{\ell}}$. We show that, up to replacing $v_{n-1}$ with another $p^{\ell}$-th  root $u_{n-1}$ of $\bar x_{n-1}$, ${\bf x}$ fulfills the parametrization of $V$. If $\bar x_{n-1}=0$, then, for all $i=1,\dots, n-2$, and for $i=n$, $F_i({\bf x})=0$ implies that $\bar y_i=0$. Hence  $u_{n-1}=0$ yields the required representation of ${\bf x}$. So suppose that $\bar x_{n-1}\ne0$. 
Let $i$ be any index with $1\leq i\leq n-2$. Condition $F_i({\bf x})=0$ implies that
$\bar y_i^{p^{\ell}}=u_i^{p^{\ell}b_i}v_{n-1}^{p^{\ell}c_i}$, i.e., 
\begin{equation}\label{yi2}
\bar y_i=u_i^{b_i}v_{n-1}^{c_i}\omega_i
\end{equation}
\noindent
for some $\omega_i\in K$ such that 
\begin{equation}\label{omegai} \omega_i^{p^{\ell}}=1.
\end{equation}
\noindent
Now let $\bar\eta$ be a primitive $p^{\ell}$-th root of 1. There is an integer $z_i$ such that $\omega_i=\bar\eta^{z_i}$. On the other hand, by assumption (I) there is an integer $w_i$ such that $c_iw_i\equiv z_i$ (mod $p^{\ell}$). Set $\eta=\bar\eta^{w_i}$. Then 
\begin{equation}\label{omegaieta}
\omega_i=\eta^{c_i}.
\end{equation}
\noindent
If $n=3$, then the only index $i$ to be considered is $i=1$, and setting $u_{n-1}=v_{n-1}\eta$ will give
$$\bar y_1=u_1^{b_1}u_{n-1}^{c_1}.$$
\noindent
Now let $n\geq 4$. We show that the same choice of parameter $u_{n-1}$ yields the required representation for all $\bar y_1, \dots, \bar y_{n-2}$. This is trivially true if $u_i=0$ (i.e., $\bar y_i=0$) for all $i=1,\dots, n-2$, or $u_i=\bar y_i=0$ for all but one of these indices. So suppose that, for two indices $i,j$, with $1\leq i<j\leq n-2$, we have $u_i\ne0$ and $u_j\ne0$, which, under our present assumption that $u_{n-1}\ne 0$, is equivalent to  $\bar y_i\ne 0$ and $\bar y_j\ne 0$. 
 Then, by (\ref{yi2}),  $G_{ij}({\bf x})=0$ implies that
\begin{equation}\label{yiyi+1} u_i^{b_ig_{ij}}v_{n-1}^{c_ig_{ij}}\omega_i^{g_{ij}}u_j^{b_jg_{ji}}v_{n-1}^{c_jg_{ji}}\omega_j^{g_{ji}}=u_i^{b_ig_{ij}}u_j^{b_jg_{ji}}v_{n-1}^{r_{ij}p^{m_{ij}}}.\end{equation}
\noindent
If we consider (\ref{gi}) and cancel equal terms on both sides of (\ref{yiyi+1}), we obtain
\begin{equation}\label{omegaiomegai+1}\omega_i^{g_{ij}}\omega_j^{g_{ji}}=1.\end{equation}
\noindent
Now (\ref{omegaieta}) and (\ref{omegaiomegai+1}) imply that 
\begin{equation}\label{cizi}\eta^{c_ig_{ij}}=\omega_j^{-g_{ji}}.\end{equation}
\noindent
On the other hand, by (\ref{gi}) we have
$${c_ig_{ij}}\equiv-c_jg_{ji}\quad(\mbox{mod }p^{\ell}),$$
\noindent
so that
$$\eta^{c_ig_{ij}}=\eta^{-c_jg_{ji}},$$
\noindent
i.e., by (\ref{cizi}),
$$\omega_j^{g_{ji}}=\eta^{c_jg_{ji}}.$$
\noindent
Since $p$ does not divide $g_{ji}$, there is an integer $q$ such that $g_{ji}q\equiv 1$ (mod $p^{\ell}$).  Therefore, by (\ref{omegai}) and (\ref{omegaieta}), applied to the index $j$,  
\begin{equation}\label{omegai+1eta}
\omega_j=\omega_j^{g_{ji}q}=\eta^{c_jg_{ji}q}=\eta^{c_j}.
\end{equation}
\noindent
Note that (\ref{omegaieta}), together with (\ref{omegai+1eta}), implies that $\omega_i=\eta^{c_i}$ holds for all indices $i=1,\dots, n-2$ such that $\bar y_i\ne0$. Set $u_n=v_n\eta$. Then, in view of (\ref{yi2}), for all these indices we have
\begin{equation}\label{yiui} \bar y_i=u_i^{b_i}u_{n-1}^{c_i}.\end{equation}
\noindent
This obviously also holds when $u_i=0$. 
It remains to show that $\bar y_n$ has the required form. This is certainly true if $\bar x_n=0$: then $u_n=0$, and $F_n({\bf x})=0$ implies that $\bar y_n=0$.\newline
 So suppose that $\bar x_n\ne 0$. 
 Condition $F_n({\bf x})=0$ implies that $\bar y_n^{p^{\ell}}=u_{n-1}^{p^{\ell}}u_n^{p^{\ell}}$, i.e., 
\begin{equation}\label{yn2} \bar y_n=u_{n-1}u_n\omega,
\end{equation}
\noindent
for some $\omega\in K$ such that 
\begin{equation}\label{omega} \omega^{p^{\ell}}=1.
\end{equation}
\noindent
If $\bar x_i=0$ (i.e., $u_i=0$) for all $i=1,\dots, n-2$, replace $u_{n-1}$ with $u_{n-1}\omega$. This will produce in (\ref{yn2}) the required representation for $\bar y_n$, and will not affect the remaining entries of ${\bf x}$. So assume that $\bar x_j\ne 0$ for some index $j$, $1\leq j\leq n-2$.   
We have that $H_j({\bf x})=0$, together with (\ref{yiui}) and (\ref{yn2}),  implies that
$$u_j^{b_jh_j}u_{n-1}^{c_jh_j}u_{n-1}^{k_j}u_n^{k_j}\omega^{k_j}=
u_j^{b_jh_j}u_{n-1}^{r_jp^{m_j}}u_n^{as_j}u_n^{dt_j}.$$
\noindent
In view of (\ref{hiki}) and (\ref{siti}), applied to the index $j$, if we cancel equal terms on both sides of the above equation, we obtain
\begin{equation}\label{omegaomegai} 
\omega^{k_j}=1,\end{equation}
\noindent
Since $p$ does not divide $k_j$, (\ref{omega}) and (\ref{omegaomegai}) imply that 
$$
\omega=1.$$
\noindent
Thus by (\ref{yn2}),
$$
\bar y_n=u_{n-1}u_n.
$$
\noindent
This  completes the proof.
\par\medskip\noindent
\section{The binomial arithmetical rank}
We have just proven that $V$ can always be set-theoretically defined by $4$ binomial equations if $n=3$ and  by $2n-2+{n-2\choose 2}$ binomial equations if $n\geq 4$. We now show that these numbers cannot be made smaller if char\,$K\ne p$.
\begin{proposition}\label{bar} If {\rm char}$\,K\ne p$, then 
$${\rm bar}\,V=\left\{\begin{array}{ll} 4&\mbox{if }n=3;\\
2n-2+{n-2\choose 2}&\mbox{if }n\geq4.\\
\end{array}\right.$$
\end{proposition}
\demo In view of the above remark, we only have to prove the inequality $\geq$. Let ${\cal B}$ be a set of binomials such that $V$ is set-theoretically defined by the vanishing of all elements of ${\cal B}$. Of course we may assume that all elements of ${\cal B}$ are irreducible. Since, by Hilbert's Nullstellensatz,  for all $i=1,\dots, n-2$, $F_i$ belongs to the radical of the ideal generated by ${\cal B}$ in $R$, one binomial of ${\cal B}$ is monic in $y_i$. By Lemma \ref{preliminary2} it follows that the binomials $F_1, \dots, F_{n-1}$ defined in (\ref{Fi}) and (\ref{Fn-1}) belong to ${\cal B}$, together with some binomial $F'_n$ monic in $y_n$. 
Let $i$ be any index with $1\leq i\leq n-2$, and let $\eta$ be a  primitive $p^{\ell}$-th root of unity.
 Consider ${\bf x}=(\bar x_1,\dots, \bar x_n,\bar y_1,\dots, \bar y_n)\in K^{2n}$, where $\bar x_i=\bar x_{n-1}=\bar x_n=1$, $\bar y_i=\eta$ and $\bar y_{n-1}=\bar y_n=1$, whereas the remaining entries are zero. Then $F_i({\bf x})=F_n({\bf x})=0$.  Suppose for a contradiction that every  $B\in {\cal B}\setminus\{F_i, F_n\}$ has a monomial $M$ such that
\begin{list}{}
\item{(a)} supp\,$(M)\subset\{x_{n-1}, x_n, y_{n-1}, y_n\}$, or
\item{(b)} supp\,$(M)\not\subset\{x_i, x_{n-1}, x_n, y_i,  y_{n-1}, y_n\}$.
\end{list}
\par\bigskip\noindent
Let $M'$ be the other monomial of such a binomial $B$. In case (a), by Lemma \ref{preliminary1} (i) it follows that supp\,$(M')\subset\{x_{n-1}, x_n, y_{n-1}, y_n\}$, so that $M({\bf x})=M'({\bf x})=1$, and, consequently, $B({\bf x})=0$. In case (b), we have that, for some index $k\ne i$ with $1\leq k\leq n-2$, $M$ is divisible by $x_k$ or $y_k$. By Lemma \ref{preliminary1} (i) the same is true for $M'$. Since $\bar x_k=\bar y_k=0$, we have that $M({\bf x})=M'({\bf x})=0$, and thus $B({\bf x})=0$. Therefore we have that, in any case, $B({\bf x})=0$ for all ${\bf x}\in{\cal B}$. We show that, however, ${\bf x}\not\in V$. If ${\bf x}\in V$, then ${\bf x}$ would fulfill the parametrization of $V$ for a suitable choice of parameters $u_1,\dots, u_n$. Then necessarily $u_i=\bar x_i$, $u_n=1$, and   $u_{n-1}=\eta^k$ for some integer $k$. Consequently, we would have
$$\eta=\bar y_i=u_i^{b_i}u_{n-1}^{c_i}=\eta^{kc_i},$$
$$1=\bar y_n=u_nu_{n-1}=\eta^k.$$
\noindent
 which implies $\eta=1$, against the definition of $\eta$. This shows that, for all $i=1,\dots, n-2$, there is a binomial $H'_i\in{\cal B}\setminus\{F_i, F_n\}$ such that none of its monomials fulfills (a) or (b), i.e., both its monomials are of the form $x_i^{\alpha_i}x_{n-1}^{\alpha_{n-1}}x_n^{\alpha_n}y_i^{\beta_i}y_{n-1}^{\beta_{n-1}}y_n^{\beta_n}$, with $\alpha_i>0$ or $\beta_i>0$. From Lemma \ref{preliminary1} and Lemma \ref{preliminary2} it follows that the remaining exponents are not all zero. This suffices to prove the first case of the claim:   since $\{F_1, F_2, F'_3, H'_1\}\subset {\cal B}$, we have
$$\mbox{for }n=3,\qquad\vert{\cal B}\vert \geq 4.$$
\noindent
Now suppose that $n\geq 4$. Note that, by construction, the $H'_i$ are $n-2$ pairwise distinct binomials. Let $i,j$ be indices such that $1\leq i<j\leq n-2$. There are  integers $d_{ij}, d_{ji}$ such that 
\begin{equation}\label{didj} c_id_{ij}-c_jd_{ji}=\gcd(c_i,c_j).
\end{equation}
\noindent
Consider ${\bf x}\in K^{2n}$ where $\bar x_i=\bar x_j=\bar x_{n-1}=1$ and $\bar y_i=\eta^{d_{ij}}$, $\bar y_j=\eta^{d_{ji}}$, whereas the remaining entries are zero. Then $F_i({\bf x})=F_j({\bf x})=0$. Suppose, for a contradiction,  that for all  $B\in{\cal B}\setminus\{F_i,F_j\}$,  the support of neither monomial of $B$ is contained in $\{x_i, x_j, y_i,  y_j\}$. Let $B$ be any such monomial, and let $M, M'$ be its monomials. Then, in view of Lemma \ref{preliminary1}, up to interchanging $M$ and $M'$, we have one of the following cases:
\begin{list}{}
\item{(a)} for some index $k$ with $1\leq k\leq n-2$, $k\ne i,j$, $M$ is divisible by $x_k$, and $M'$ is divisible by $y_k$;
\item{(b)}  $M$ is divisible by one of the indeterminates $x_n, y_{n-1}, y_n$, and $M'$ is divisible by one of the remaining two;
\item{(c)} $M$ is divisible by $x_{n-1}$; in this case, by Lemma \ref{preliminary1} (iii),  $M'$ is divisible by one of the indeterminates $y_1,\dots, y_{n-2}, y_n$; since supp\,$(M')\not\subset\{y_i, y_j\}$, this takes us back to case (a) or (b).
\end{list}
\par\smallskip\noindent
In all the above cases $M({\bf x})=M'({\bf x})=0$. We conclude that $B({\bf x})=0$ for all $B\in{\cal B}$.  Once again we show that assuming  ${\bf x}\in V$ leads to a contradiction.  In fact, under this assumption, ${\bf x}$ fulfills the parametrization of $V$ with $u_i=u_j=1$ and with  $u_{n-1}=\eta^k$ for some integer $k$. Consequently, 
$$ u_i^{b_i}u_{n-1}^{c_i}=\eta^{kc_i},\qquad\qquad
u_j^{b_j}u_{n-1}^{c_j}=\eta^{kc_j},$$
\noindent
Hence $\bar y_i=u_i^{b_i}u_{n-1}^{c_i}$ and $\bar y_j=u_j^{b_j}u_{n-1}^{c_j}$ imply 
$$
\eta^{d_{ij}}=\eta^{kc_i},\qquad\qquad
\eta^{d_{ji}}=\eta^{kc_j},$$
\noindent
which is equivalent to
$$
d_{ij}\equiv kc_i,\mbox{ (mod $p^{\ell}$)},\qquad\qquad
d_{ji}\equiv kc_j\mbox{ (mod $p^{\ell}$)}.$$
\noindent
Therefore
$$c_jd_{ij}\equiv c_id_{ji}\mbox{ (mod $p^{\ell}$)},$$
\noindent
which is incompatible with (\ref{didj}), since, in view of assumption (I), $p$ does not divide $\gcd(c_i,c_j)$.
This  shows that ${\bf x}\not\in V$ and provides the required contradiction. We conclude that ${\cal B}$ must contain, for all indices $i,j$ such that $1\leq i<j\leq n-2$, a binomial $G'_{ij}$, other than $F_i,F_j$, such that  one of its monomials is of the form
$x_i^{\gamma_i}x_j^{\gamma_j}y_i^{\delta_i}y_j^{\delta_j}$. By Lemma \ref{preliminary1} (i), Lemma \ref{preliminary2} and irreducibility, it follows that one of $\gamma_i, \delta_i$ and one of $\gamma_j, \delta_j$ are positive. It follows that the $G'_{ij}$ are ${n-2\choose 2}$ pairwise distinct binomials. It also evident that the sets $\{F_1,\dots, F_{n-1}, F'_n\}$, $\{G_{12},\dots, G_{n-2\,n-3}\}$ and $\{H'_1,\dots, H'_{n-2}\}$ are pairwise disjoint.  Therefore 
$$\mbox{for }n\geq 4,\qquad\vert{\cal B}\vert \geq n+{n-2\choose 2}+n-2= 2n-2+{n-2\choose 2}$$
\noindent
This completes the proof. 

\section{A lower bound for the arithmetical rank}
In this section we give a lower bound for ara\,$V$ when char$\,K\ne p$.  
  We will use the following result, which is due to Newstead and is quoted from \cite{BS}. It is based on \'etale cohomology ($\het$). We refer to \cite{M} or to \cite{M2} for the basic notions.
\begin{lemma}\label{Newstead} Let
$W\subset\tilde W$ be affine varieties. Let $d=\dim\tilde
W\setminus W$. If there are $s$ equations $F_1,\dots, F_s$ such
that $W=\tilde W\cap V(F_1,\dots,F_s)$, then 
$$\het^{d+i}(\tilde W\setminus W,{\bz}/r{\bz})=0\quad\mbox{ for all
}i\geq s$$ and for all $r\in{\bz}$ which are prime to {\rm char}\,$K$.
\end{lemma}
\par\medskip\noindent
We prove the following result. 
\begin{proposition}\label{ara} If {\rm char}\,$K\ne p$, then {\rm ara}\,$V\geq 2n-2$.
\end{proposition}
\demo  We have to show that $V$ cannot be defined set-theoretically by $2n-3$ equations. According to Lemma \ref{Newstead} this is true if 
$$\het^{4n-3}(K^{2n}\setminus V, \bz/p\bz)\neq0.$$
\noindent
Since $K^{2n}\setminus V$ is non singular, by Poincar\'e Duality (see \cite{M}, Corollary 11.2, p.~276), this is equivalent to
\begin{equation}\label{claim}\hc^3(K^{2n}\setminus V, \bz/p\bz)\neq0,\end{equation}
\noindent where $\hc$ denotes cohomology with compact support.
For the sake of simplicity, we shall henceforth omit the coefficient group $\bz/p\bz$. According to \cite{M}, Remark 1.30, p.~94, There is an exact sequence: 
\begin{equation}\label{exact1}\hc^2(K^{2n})\rightarrow\hc^2(V)\rightarrow\hc^3(K^{2n}\setminus V)\rightarrow\hc^3(K^{2n}).\end{equation}
\noindent
Recall that 
\begin{equation}\label{K}\hc^{i}(K^m)\simeq\left\{\begin{array}{ll}{\bz}/p{\bz}&\mbox{ for }i=2m\\
0&\mbox{ otherwise }\end{array}\right.\end{equation}
\noindent
See \cite{M2}, Example 16.3, pp.~98--99, together with Poincar\'e Duality, for a proof.
Hence the left and the right group in (\ref{exact1}) are zero, so that the middle map is a group isomorphism. Thus our claim (\ref{claim}) is equivalent to 
\begin{equation}\label{V} \hc^2(V)\ne0.\end{equation}
\noindent
\par\smallskip\noindent
The proof of (\ref{V}) needs some preparation.
Consider the following morphism of schemes:
$$\varphi: K^n\longrightarrow  V$$
$$(u_1,\dots, u_n)\mapsto(u_1, \dots, u_{n-2}, u_{n-1}^{p^{\ell}}, u_n^a, u_1^{b_1}u_{n-1}^{c_1},\dots, u_{n-2}^{b_{n-2}}u_{n-1}^{c_{n-2}}, u_n^d, u_{n-1}u_n).$$
\noindent
Let $W$ be the subvariety of $K^n$ defined by $u_1u_{n-1}=\cdots=u_{n-2}u_{n-1}=u_n=0$. We show that  $\varphi$ induces by restriction  
an isomorphism of schemes:
$$\tilde\varphi: K^n\setminus W\longrightarrow V\setminus\varphi(W).$$
\noindent
Note that $V\setminus\varphi(W)$ is the union of the open subsets
$$V_i=\{{\bf x}\in V\vert y_i\ne 0\}\qquad(i=1,\dots, n-1).$$
\noindent
Moreover, for  all ${\bf x}\in V$ and all $i=1,\dots, n-2$, 
$$y_i\ne 0\quad\mbox{is equivalent to}\quad x_i\ne 0\mbox { and }x_{n-1}\ne0,$$
\noindent
and
$$y_{n-1}\ne0\quad\mbox{is equivalent to}\quad x_n\ne0.$$
\noindent
Thus we have
$$ U_i=\varphi^{-1}(V_i)=\{{\bf u}\in K^n\vert u_iu_{n-1}\ne 0\}\qquad(i=1,\dots, n-2),$$
\noindent
and
$$U_{n-1}=\varphi^{-1}(V_{n-1})=\{{\bf u}\in K^n\vert u_n\ne 0\}.$$
\noindent
By assumption (II) of Section 1 there are integers $r,s$ such that $ar+ds=1$; by assumption (I) there are, for all $i=1,\dots, n-2$, integers $v_i, w_i$ such that $c_iv_i+p^{\ell}w_i=1$.
The following morphisms of schemes are inverse to the restrictions of $\varphi$ to $U_i$.
$$\qquad\qquad\qquad V_i\longrightarrow U_i\qquad(i=1,\dots, n-2)$$
$$(x_1,\dots, x_n, y_1,\dots, y_n)\mapsto \left(x_1,\dots, x_{n-2}, \frac{y_i^{v_i}x_{n-1}^{w_i}}{x_i^{b_iv_i}}, \frac{y_nx_i^{b_iv_i}}{y_i^{v_i}x_{n-1}^{w_i}}\right),$$
\hskip.1truecm
$$\!\!\!\!\!\!\!\!\!\!\!\!\!\!\!\!\!\!\!\!\!\!\!\!V_{n-1}\longrightarrow U_{n-1}$$
$$(x_1,\dots, x_n, y_1,\dots, y_n)\mapsto \left(x_1,\dots, x_{n-2}, \frac{y_n}{x_n^ry_{n-1}^s}, x_n^ry_{n-1}^s\right).$$
\noindent
We have just proven that $\varphi$ is an isomorphism. Hence, for all indices $i$, it induces an isomorphism in cohomology with compact support:
\begin{equation}\label{iso}\varphi_i:\hc^i(V\setminus\varphi(W))\mathop{\longrightarrow}^{\simeq}\hc^i(K^n\setminus W).\end{equation}
\par\smallskip\noindent
Now consider the subvariety $Y$ of $W$ defined by $u_1=\cdots=u_{n-2}=u_n=0$. Then $Y$ can be identified with $K$, and
$W\setminus Y$ with the set of all points in $K^{n-1}$ such that $u_{n-1}=0$, whereas not all $u_i$ with $1\leq i\leq n-2$ are zero; in other words,  $W\setminus Y$ can be identified with $K^{n-2}\setminus\{0\}$. It can be easily seen that also $\varphi(Y)$ ( which is a closed subset of $\varphi(W)$) and $\varphi(W)\setminus\varphi(Y)$ can be identified with $K$ and $K^{n-2}\setminus\{0\}$ respectively.  From (\ref{K}) and the long exact sequence in \cite{M}, Remark 1.30, p.~94, it easily follows that
\begin{equation}\label{Kstar}\hc^{i}(K^m\setminus\{0\})\simeq\left\{\begin{array}{ll}{\bz}/p{\bz}&\mbox{ for }i=1,2m\\
0&\mbox{ otherwise }\end{array}\right..\end{equation}
\noindent
Now, in view of the above identifications, (\ref{K}) and (\ref{Kstar}) imply
\begin{equation}\label{Y}\hc^i(\varphi(Y))\simeq\hc^{i}(Y)\simeq\left\{\begin{array}{ll}{\bz}/p{\bz}&\mbox{ for }i=2\\
0&\mbox{ otherwise }\end{array}\right.,\end{equation}
\noindent
and
\begin{equation}\label{WY}\hc^{i}(\varphi(W)\setminus\varphi(Y))\simeq\hc^{i}(W\setminus Y)\simeq\left\{\begin{array}{ll}{\bz}/p{\bz}&\mbox{ for }i=1, 2n-4\\
0&\mbox{ otherwise }\end{array}\right..\end{equation}
\noindent
In the sequel we shall use the following exact sequences:
\begin{equation}\label{exact2}\hc^{i}(K^n)\rightarrow\hc^{i}(W)\rightarrow\hc^{i+1}(K^n\setminus W)\rightarrow\hc^{i+1}(K^n),\end{equation}
\noindent
\begin{equation}\label{exact3}\hc^{i-1}(Y)\rightarrow\hc^{i}(W\setminus Y)\rightarrow\hc^{i}(W)\rightarrow\hc^{i}(Y)\rightarrow\hc^{i+1}(W\setminus Y).\end{equation}
\noindent
There is a similar sequence obtained from (\ref{exact3}) by replacing $W$ and $Y$ with $\varphi(W)$ and $\varphi(Y)$ respectively. 
We are now ready to prove claim (\ref{V}). 
Consider sequence (\ref{exact3}) for $i=2$. We have
\begin{equation}\label{exact4}\begin{array}{ccccccccc}
\hc^1(Y)&\rightarrow&\hc^2(W\setminus Y)&\rightarrow&\hc^2(W)&\rightarrow&\hc^2(Y)&\rightarrow&\hc^3(W\setminus Y)\\
\|&&&&&&|\wr&&\|\\
0&&&&&&\bz/p\bz&&0
\end{array},\end{equation}
\noindent
where the equalities and the isomorphisms follow from (\ref{Y}) and (\ref{WY}). From (\ref{WY}) we also have that
$$\hc^2(W\setminus Y)\simeq\left\{\begin{array}{ll}{\bz}/p{\bz}&\mbox{ if }n=3\\
0&\mbox{ if }n\geq4\end{array}\right..$$
\noindent
In view of (\ref{exact4}) we deduce that 
$$\vert\hc^2(W)\vert=\left\{\begin{array}{ll}p^2&\mbox{ if }n=3\\
p&\mbox{ if }n\geq4\end{array}\right..$$
\noindent
A similar result holds for $\hc^2(\varphi(W))$.
Consequently,
\begin{equation}\label{cardinality} \vert\hc^2(\varphi(W))\vert=\vert\hc^2(W)\vert\ne0.\end{equation}
\noindent
On the other hand, from the exact sequence (\ref{exact2}), for $i=2$, we have
$$\begin{array}{ccccccc}
\hc^2(K^n)&\rightarrow&\hc^2(W)&\rightarrow&\hc^3(K^n\setminus W)&\rightarrow&\hc^3(K^n)\\
\|&&&&&&\|\\
0&&&&&&0
\end{array},
$$
where the equalities are a consequence of (\ref{K}). We thus have an isomorphism:
\begin{equation}\label{iso2}\hc^2(W)\mathop{\longrightarrow}^{\simeq}\hc^3(K^n\setminus W).\end{equation}
\noindent
Next we show that the following map, induced in cohomology by the restriction of $\varphi$ to $W$,
$$\varphi'_2: \hc^2(\varphi(W))\longrightarrow \hc^2(W)$$
\noindent
is not injective. 
It is well known that the restriction of $\varphi$ to $Y$
$$\varphi_|:Y\longrightarrow \varphi(Y)$$
$$u_{n-1}\mapsto u_{n-1}^{p^{\ell}}$$
\noindent
induces multiplication by $p^{\ell}$ in cohomology (see \cite{M2}, Remark 24.2 (f), p.~135). Thus $\varphi$ gives rise to the following commutative diagram
$$\begin{array}{ccccc}
&&\bz/p\bz&&\\
&&|\wr&&\\
\hc^2(W)&\displaystyle\mathop{\rightarrow}^{\alpha}&\hc^2(Y)&\rightarrow&0\\\\
\varphi'_2\mbox{\LARGE$\uparrow$}&&\cdot p^{\ell}\mbox{\LARGE$\uparrow$}&\\\\
\hc^2(\varphi(W))&\rightarrow&\hc^2(\varphi(Y))&\rightarrow&0\\
&&|\wr&&\\
&&\bz/p\bz&&
\end{array},$$
\noindent
where the first row is part of the exact sequence (\ref{exact4}) and the second row is similarly derived from (\ref{exact3}).  Note that multiplication by $p^{\ell}$ is the zero map. If $\varphi'_2$ were injective, in view of (\ref{cardinality}) it would also be surjective. But then so would be $\alpha\varphi'_2$; this map, however, because of commutativity, is the zero map, which is a contradiction.\par\smallskip\noindent
Finally consider the following commutative diagram
 $$\begin{array}{cccccccc}
&&\hc^2(W)&\displaystyle\mathop{\rightarrow}^{\simeq}&\hc^3(K^n\setminus W)&&\\\\
&&\varphi'_2\mbox{\LARGE$\uparrow$}&&\varphi_3\mbox{\LARGE$\uparrow$}|\wr&&&\\\\
\hc^2(V)&\rightarrow&\hc^2(\varphi(W))&\displaystyle\mathop{\rightarrow}^{\beta}&\hc^3(V\setminus\varphi(W))&&
\end{array},$$
\noindent
where the isomorphisms are those given in (\ref{iso}) and (\ref{iso2}). 
Since $\varphi'_2$ is not injective, nor is $\beta$. It follows that $\hc^2(V)\ne 0$, i.e., (\ref{V}) is true. This completes the proof.\par\medskip\noindent
In the special case where $n=3$ Proposition \ref{ara} yields ara\,$V\geq 4$; on the other hand,  by Proposition \ref{bar}, we have  bar\,$V=4$. Thus ara\,$V=$\,bar\,$V=4$. 
Moreover, according to a classical theorem, proven by Eisenbud and Evans \cite{EE}, and, independently, by Storch \cite{S}, every variety in the $N$-dimensional affine space can be defined by a system of $N$ equations. Thus the above results can be summarized as follows:
\begin{corollary}\label{corollary3} It holds:
\begin{list}{}{}
\item{(i)} {\rm ara}\,$V = n$, if {\rm char}$\,K=p$;
\item{(ii)} $2n-2\leq{\rm ara}\,V \leq 2n$,  if {\rm char}\,$K\ne p.$
\end{list}
\noindent
In (i), {\rm ara}\,$V$ defining equations can be chosen to be binomial, whereas this is possible in (ii) if and only if $n=3$, and in this case {\rm ara}\,$V=4$. 
\end{corollary}
The problem of determining ara\,$V$  when char\,$K\ne p$  and $n\geq 4$ remains open.

\end{document}